\documentclass[12pt]{article}
 \usepackage{amsmath}

\usepackage{amssymb,amsthm}
 \usepackage{mathrsfs,amsfonts}
 \usepackage{color}
 \usepackage{graphicx}
\usepackage{bbm}
\usepackage{ulem}

\newcommand{\ddr}{\mathrm{d}}
 \newtheorem{theorem}{Theorem}[section]
 \newtheorem{lemma}[theorem]{Lemma}
 \newtheorem{corol}[theorem]{Corollary}
 \newtheorem{prop}[theorem]{Proposition}
 \newtheorem{example}[theorem]{Example}
 \newtheorem{Definition}[theorem]{Definition}
 \newtheorem{con1}[theorem]{Condition}
\theoremstyle{remark}
 \newtheorem{remark}[theorem]{Remark}

 \def\blemma{\begin{lemma}\sl{}\def\elemma{\end{lemma}}}

 \def\beqlb{\begin{eqnarray}}\def\eeqlb{\end{eqnarray}}
 \def\beqnn{\begin{eqnarray*}}\def\eeqnn{\end{eqnarray*}}

 \def\d{{\rm d}}

 \def\proof{\noindent{\it Proof.~~}}\def\qed{\hfill$\Box$\medskip}

 \def\<{\langle}\def\>{\rangle}

 \def\mcr{\mathscr}\def\mbb{\mathbb}

 \def\ar{\!\!&}

\providecommand{\href}[2]{#2}
\begin{document}

%(\today)
\bigskip\bigskip
%\centerline{\Large\bf Speeds of coming down from infinity}
%\bigskip
%\centerline{\Large\bf for continuous-state nonlinear branching processes}

\centerline{\Large\bf On the entrance at infinity of Feller processes}
%\bigskip
\centerline{\Large\bf with no negative jumps}

\bigskip\bigskip

\bigskip

\centerline{\large  Cl\'ement Foucart \footnote{Universit\'e Sorbonne Paris Nord and Paris 8, Laboratoire Analyse, G\'eom\'etrie $\&$ Applications, UMR 7539. Institut Galil\'ee, 99 avenue J.B. Cl\'ement, 93430 Villetaneuse, France, foucart@math.univ-paris13.fr}, Pei-Sen Li \footnote{Institute for Mathematical Sciences, Renmin University of China, 100872 Beijing, P. R. China, peisenli@ruc.edu.cn} and Xiaowen Zhou \footnote{Department of Mathematics and Statistics, Concordia University, 1455 De Maisonneuve Blvd. W., Montreal, Canada, xiaowen.zhou@concordia.ca}}

%\begin{aug}
%\author{\fnms{Cl\'ement} \snm{Foucart}\thanksref{a}\ead[label=e1]{foucart@math.univ-paris13.fr}}
%\author{\fnms{Pei-Sen} \snm{Li}\thanksref{b}\ead[label=e2]{peisenli@ruc.edu.cn}}
%\and
%\author{\fnms{Xiaowen} \snm{Zhou}\thanksref{c}\ead[label=e3]{xiaowen.zhou@concordia.ca}}
%
%\address[a]{Universit\'e Paris 13, Laboratoire Analyse, G\'eom\'etrie $\&$ Applications, UMR 7539 Institut Galil\'ee, 99 avenue J.B. Cl\'ement, 93430 Villetaneuse, France.
%\printead{e1}}
%
%\address[b]{Institute for Mathematical Sciences, Renmin University of China, 100872 Beijing, P. R. China.
%\printead{e2}}
%
%\address[c]{Department of Mathematics and Statistics, Concordia University, 1455 De Maisonneuve Blvd. W., Montreal, Canada
%\printead{e3}}
%
%
%\runauthor{C. Foucart et al.}
%
%\affiliation{Some University and Another University}

%\end{aug}

%%%%%%%%%%%%%%%%%%%%%%%%%%%%%%%%%%%%%%%

\begin{abstract}
Consider a non-explosive positive Feller process with no negative jumps. It is shown in this note that when infinity is an entrance boundary, in the sense that the entrance times of the process remain bounded when the initial value tends to infinity, the process admits a Feller extension on the compactified state space $[0,\infty]$. Moreover, when started from infinity, the extended Markov process on $[0,\infty]$ leaves infinity instantaneously and stays finite,  almost-surely. Arguments are adapted from a proof given by O. Kallenberg \cite{Kallenberg} for diffusions. We also show that the process started from $x$ converges weakly towards that started from infinity in the Skorokhod space, when $x$ goes to infinity.
\end{abstract}
\section{Introduction}
In the last decade, a regain of attention in the literature has been paid to the study of one-dimensional Markov processes at their boundaries.  Such studies are of interest for instance when the process is representing the size of a random population. There is a large body of literature on this topic and we refer the reader for instance to the recent articles of Bansaye et al. \cite{2017arXiv171108603B}, \cite{MR3595771}, D\"oring and Kyprianou \cite{2018arXiv180201672D}, Foucart \cite{Foucart},  Le and Pardoux \cite{LePardoux15} and \cite{LePardoux19}, Li \cite{2016arXiv160909593L}  and Li et al. \cite{2017arXiv170801560L}. See also the monography of Pardoux \cite[Chapters 6 and 8]{Pardouxbook}.

In the latter cited works, the authors share the same interest to understand how the population behaves when the initial size takes arbitrarily large values. Several behaviors are possible. When the dynamics prevents the process to explode  (i.e. to hit $\infty$ in finite time) and yet allows the process to ``start from infinity"; the boundary is called an \textit{entrance}. Such a phenomenon is ubiquitous in population models with self-regulation properties or in models of statistical  physics.

Several definitions of entrance boundary are given in the literature. The process is often said to \textit{come down from infinity}, when the first entrance times of compact sets are, in some sense, uniformly bounded with respect to the initial state of the process.

A rigorous definition in this vein can be found in Revuz and Yor \cite[Chapter 3]{RevuzYor} and Kallenberg \cite[Chapter 23]{Kallenberg} and is stated as follows.
\begin{Definition}
Let $(X_t,t\geq 0)$ be a positive Markov process. 
%For any $x\in [0,\infty)$, denote by $\mathbb{P}_x$, its law started from $x$. 
The boundary $\infty$ is said to be an \textit{instantaneous entrance boundary} for the process $(X_t,t\geq 0)$ if the process does not explode and
\begin{equation}\label{0.5'} \forall t>0, \underset{b\rightarrow \infty}{\lim}\underset{x\rightarrow \infty}{\liminf}\ \mathbb{P}_x(T_b\leq t)=1
\end{equation}
where $\mbb{P}_x(\cdot):=\mbb{P}(\cdot |X_0= x)$ and $T_b:=\inf\{t\geq 0; X_t< b\}$ for all $b\in \mathbb{R}_+$, with the convention $\inf\emptyset:=\infty$.
\end{Definition}
Although the meaning of condition \eqref{0.5'} is intuitively clear, it does not guarantee a priori, that  the process $X$ can be started at $\infty$.

Kallenberg \cite{Kallenberg} has designed elegant arguments for diffusions ensuring that indeed, if the assumption \eqref{0.5'} holds, then the diffusion process $(X_t,t\geq 0)$ has a regular version started from infinity. The main purpose of this note is to observe that, more generally if the process under study has no negative jumps and satisfies a Feller property on $E:=[0,\infty)$, then the assumption  \eqref{0.5'} ensures that its semigroup can be extended into a Feller semigroup on $\bar{E}=[0,\infty]$. As a direct consequence, any Feller process $(X_t,t\geq 0)$ on $E$, with no negative jumps, satisfying \eqref{0.5'} admits a c$\grave{\text{a}}$dl$\grave{\text{a}}$g extension  started from $\infty$.

%The arguments to prove this result are adapted from those of Kallenberg in \cite[Theorem 23.13]{Kallenberg}. We will also establish an important regularity result for the first passage time $T_b$ and the convergence in law of the processes started from $x$ to that started from infinity, when $x$ goes to $\infty$, in the Skorokhod sense.
%\\

The following equivalent conditions for (\ref{0.5'}) to hold, have already appeared in the literature, see e.g. \cite[Proposition 2.12]{2017arXiv170801560L} and \cite[Theorem 1.11]{2016arXiv160909593L}. The proof is deferred to Section \ref{prooflemmaequivcdi}. Denote by $\mbb{E}_x$ the corresponding expectation under $\mbb{P}_x.$

\begin{lemma}\label{equivcdi} \
Consider a positive strong Markov process $(X_t,t\geq 0)$ with no negative jumps. For any $b\in E$, recall $T_b$ its first passage time below $b$. The following statements are equivalent:
\begin{itemize}
\item[(a)] Condition \eqref{0.5'} holds.
\item[(b)] For all large enough $b$, $\sup_{x\geq b}\mathbb{E}_x(T_b)<\infty$.
\item[(c)] $\underset{b\rightarrow \infty}{\lim}\ \underset{x\rightarrow \infty}{ \lim }\mathbb{E}_x(T_b)=0$.
\end{itemize}
\end{lemma}
\section{Main results }\label{Fellerextension}
Recall $E:=[0,\infty)$ and $\bar{E}:=[0,\infty]$. Denote by $C_b(E)$ the space of continuous bounded functions on $E$ and by $C(\bar{E})$ the space of continuous functions on $\bar{E}$, namely those in $C_b(E)$ with finite limits at $\infty$.

Consider a Markov process $(X_t,t\geq 0)$ with state space $E$. Denote by $(P_t)_{t\geq 0}$ its semigroup. Assume that
\begin{itemize}
%\item[\textit{(i})] the process has no negative jumps,
\item[\textit{(i)}] for any $f\in C_b(E)$, $P_tf\in C_b(E)$,
\item[\textit{(ii)}] for any $f\in C_b(E)$ and $x\in E$, $P_tf(x)\underset{t\rightarrow 0}{\longrightarrow} f(x)$.
\end{itemize}
We shall say in the sequel that a process is Feller on the space $\bar{E}$, if its semigroup satisfies \textit{(i)} and \textit{(ii)} for any function $f\in C(\bar{E})$ and any $x\in \bar{E}$.
%$X$ can be extended into a Feller process over $[0,\infty]$ if for any $f$
%We now provide some fundamental results for general Feller Markov processes with no negative jumps.
\begin{remark}
%\begin{itemize}
%\item[a)]
Several definitions of ``Feller processes" coexist in the literature. We are following the definition given in Kallenberg's book \cite[Chapter 19, p.369]{Kallenberg}. In order to deal with an entrance boundary at $\infty$, the usual class $C_0$ of continuous functions vanishing at $\infty$ cannot be used. Recall however that when the state space is compact, as  the case for $\bar{E}$, $C_0$ can be replaced by $C(\bar{E})$. Moreover, according to \cite[Theorem 19.6]{Kallenberg}, when the state-space is compact, the strong continuity at $0$ of the semigroup is equivalent to the pointwise one.
%\item[b)]
\end{remark}
%We mention that the state-space $E=[0,\infty)$ does not play any specific role and can be any real interval $(\ell, r)$ by replacing $\infty$ by $r$ everywhere in Definition \ref{0.5'} and in the forthcoming Theorem \ref{extension}, Proposition \ref{momentunderPinfty}, Theorem \ref{weakconv} and Corollary \ref{cor}.
%\end{itemize}

\begin{theorem}\label{extension}
%Consider a {\blue{Markov process}} $(X_t,t\geq 0)$ {\blue on $E$ with no negative jumps} satisfying (i) and (ii). If the process does not explode and the entrance boundary condition \eqref{0.5'} is satisfied, 
Consider a non-explosive Markov process $(X_t,t\geq 0)$ on $E$, with no negative jumps, satisfying (i), (ii) and  the entrance boundary condition \eqref{0.5'}. Then $(X_t,t\geq 0)$ can be extended into a Feller process over $\bar{E}$ such that under $\mathbb{P}_\infty$, it starts from $\infty$, leaves it instantaneously and stays finite, almost-surely:
\vspace*{-3mm}
\begin{center}
$\mathbb{P}_\infty(X_0=\infty)=1$ and $\mathbb{P}_\infty(\forall t>0, X_t<\infty)=1$.
\end{center}
\end{theorem}
\begin{remark}
Theorem \ref{extension} holds true also for a process $(X_t,t\geq 0)$ valued in $\mathbb{N}$ with the skip-free property (i.e. negative jumps of size at most $-1$). It is also worth mentioning that the state-space $E=[0,\infty)$ does not play any specific role here and can be any real interval $(\ell, r)$ by replacing $\infty$ by the right-end boundary $r$ everywhere in Definition \ref{0.5'}, in Theorem \ref{extension},  and in the forthcoming Proposition \ref{momentunderPinfty}, Theorem \ref{weakconv} and Corollary \ref{cor}.
\end{remark}
%\end{remark}
%Theorem \ref{extension} applies for instance to
%{\blue{diffusions, see e.g. \cite[Chapter 23]{Kallenberg} where an explicit condition equivalent to \eqref{0.5'} is given in terms of the scale function and the speed measure of the diffusion. We refer also to Foucart et al. \cite[Theorem 3.1]{timechangedsplp} where time-changed spectrally positive L\'evy processes started from infinity are studied and an explicit  condition equivalent to \eqref{0.5'} is found. Other examples will be discussed in the sequel.}}
Theorem \ref{extension} applies for instance to
diffusions, see e.g. \cite[Chapter 23]{Kallenberg} where an explicit condition equivalent to \eqref{0.5'} is given in terms of the scale function and the speed measure of the diffusion process. It is also invoked in Foucart et al. \cite[Section 2.1]{timechangedsplp} for defining certain time-changed spectrally positive L\'evy processes, started from infinity. An explicit condition equivalent to \eqref{0.5'} is found for this class of processes in \cite[Theorem 3.1]{timechangedsplp}. Other examples will be discussed in the sequel, see Corollary \ref{NLCSBPstartedfrominfinity}.

The next proposition states a regularity property of the first entrance times $T_b$ under the probability measures $(\mathbb{P}_x,x\geq 0)$. Denote by $\mathbb{E}_\infty$ the expectation under $\mathbb{P}_\infty$.

\begin{prop}\label{momentunderPinfty} Suppose that the assumptions of Theorem \ref{extension} hold. Let $h$ be a continuous function on $[0,\infty)$ that is either bounded or nonnegative and increasing. Then
\begin{itemize}
\item[(a)]
for any $\theta>0$, there exists $b_\theta>0$, such that for all $b\geq b_\theta$,
%\begin{equation*}
$\mathbb{E}_\infty(e^{\theta T_b})<\infty;$
%\end{equation*}
 \item[(b)] for any $b>0$, if $T_b<\infty$ $\mathbb{P}_\infty$-a.s. and $\mathbb{E}_\infty(h(T_b))<\infty$, then
 \begin{equation}\label{convergencemoment} \mathbb{E}_x (h(T_b))\underset{x\rightarrow \infty}{\longrightarrow} \mathbb{E}_\infty (h(T_b)).
\end{equation}
\end{itemize}
\end{prop}
\noindent In particular, we see from Proposition \ref{momentunderPinfty}-(a) that, under $\mathbb{P}_\infty$, one can find a large enough $b$ such that $T_b$ has  moments of all orders. The convergence \eqref{convergencemoment} holds for instance with the function $h(x)=x^{n}$ for any $n\in \mathbb{N}$. Such a convergence of moments is crucial when one wants to study the process started from infinity, from the sample paths started from large but finite levels. We refer to \cite{2017arXiv171108603B}, \cite{MR3595771}  and \cite{timechangedsplp} where it is used for studying the speed of coming down from infinity of different Feller processes with no negative jumps.

The last theorem states a general result for Feller processes valued in a compact state space and clarifies in which sense the laws $(\mathbb{P}_x)_{x\geq 0}$ converge as $x$ goes to $\infty$ towards $\mathbb{P}_\infty$.

\begin{theorem}\label{weakconv} Assume that $X$ is a Feller process on $[0,\infty]$. Let $X^{(x)}$ be the Markov process started from $x\in [0,\infty]$ with c$\grave{\text{a}}$dl$\grave{\text{a}}$g   sample paths.  Then the family of processes $(X^{(x)})_{x\in [0,\infty)}$ converges weakly, in the Skorokhod topology, as $x\to\infty$ towards $X^{(\infty)}$.
\end{theorem}
\begin{remark} Theorem \ref{weakconv} does not require the assumption of absence of negative jumps.
\end{remark}
A direct consequence of Theorem \ref{extension} and Theorem \ref{weakconv} is the following convergence in law of the process started from $x$ towards that started from $\infty$, when $\infty$ is an entrance boundary.
\begin{corol}\label{cor}
%Assume that the {\blue entrance boundary} condition \eqref{0.5'} and the assumptions (\textit{i-iii}) are fulfilled. 
Suppose that the assumptions of Theorem \ref{extension} hold. Denote by $(X^{(x)})_{x\in [0,\infty)}$ the family of processes started from $x$. Then, $(X^{(x)})_{x\in [0,\infty)}$ converges in law, in the Skorokhod topology, as $x\to\infty$ towards $X^{(\infty)}$.
\end{corol}
Such a weak convergence was previously shown for different Markov processes, satisfying the conditions $\textit{(i)}$ and $\textit{(ii)}$, with $\infty$ as entrance boundary. See for instance Donnelly \cite{Donnelly}, Bansaye \cite[Proposition 4.4]{2015arXiv151107396B} and Bansaye et al. \cite{2017arXiv171108603B}, \cite[Lemma 2.1]{MR3595771}. The arguments for establishing the weak convergence in these latter works were requiring the stochastic monotonicity of the process in its initial values.
%This assumption is somehow replaced in Corollary \ref{cor} by the assumption of  absence of negative jumps.

%\section{Applications}
We now show how our main results apply to the so-called \textit{continuous-state nonlinear branching processes}.  Those processes were defined and studied in \cite{2017arXiv170801560L}.
%In this section, we will apply our main results to the so-called general continuous-state nonlinear branching process, which will be defined below.
Let $\{B_t\}_{t\ge0}$ be an $(\mcr{F}_t)$-Brownian motion. Write $\nu$ for a $\sigma$-finite nonnegative measure on $(0,\infty)$ such that $\int^\infty_0 (z\wedge
z^2) \,\nu(\d z)<\infty.$ Let
$\{N(\d s,\d z,\d u): s,z,u>0\}$ be an independent
$(\mcr{F}_t)$-Poisson random measure on $(0,\infty)^3$ with
intensities $\d s \,\nu(\d z)\,\d u$, and $\{\tilde N(\d s,\d z,\d
u): s,z,u>0\}$ be the corresponding compensated measure, namely $\tilde N(\d s,\d z,\d
u):= N(\d s,\d z,\d u)-\d s \,\nu(\d z)\,\d u.$ Consider the following stochastic differential equation (SDE):
\begin{equation}\label{sdeB}\begin{split} X_t^{(x)} = x&+\int_0^t\gamma_0(X_s^{(x)})\,\d
s+\int_0^t \sqrt{\gamma_1(X_{s}^{(x)})} \,\d
B_s\\
&+\int_0^t\int_{0}^\infty\int_0^{\gamma_2(X_{s-}^{(x)})} z\,\tilde{N}(\d
s, \d z, \d u). \end{split}\end{equation}
Any solution to \eqref{sdeB}, whenever it exists, is called a general continuous-state nonlinear branching process. In the two next results, we shall make the following general assumption:

\noindent Condition $\mathbb{H}$: there is a pathwise unique nonnegative  c\`adl\`ag solution $(X_t^{(x)},t\geq 0)$ to \eqref{sdeB} which is non-explosive and satisfies for any $x\leq y$, \begin{equation}\label{comparisonproperty} X_t^{(x)}\leq X_t^{(y)} \text{ for all } t\geq 0 \text{ a.s. }
\end{equation}
Existence and pathwise uniqueness are guaranteed for instance if the functions $\gamma_i$, for $i\in \{0,1,2\}$, are locally Lipschitz. The comparison property \eqref{comparisonproperty} is fulfilled as soon as $\gamma_2$ is nondecreasing. We refer to \cite[Theorem 3.1]{2017arXiv170801560L} and \cite[Proposition 2.3, Theorem 5.3 and 5.5]{FL} for more details. Sufficient conditions for non-explosion of continuous-state nonlinear branching processes are provided in \cite[Theorem 2.8]{2017arXiv170801560L}.

Any c\`adl\`ag process solution to \eqref{sdeB} has  no negative jumps and satisfies the condition \textit{(ii)}. In the following proposition, we provide a  sufficient condition on the drift $\gamma_0$ entailing the Feller property \textit{(i)}.
\begin{prop}\label{Fellerprop} Assume that Condition $\mathbb{H}$ holds. If there exists $\theta>0$ such that for any $0\leq x\leq y$, \begin{equation}\label{onesidedLipschitz}
\gamma_0(y)-\gamma_0(x)\leq \theta (y-x),
\end{equation}
then the semigroup of $(X_t,t\geq 0)$ satisfies the Feller property \textit{(i)}.
\end{prop}

Combining Condition $\mathbb{H}$, the one-sided Lipschitz condition \eqref{onesidedLipschitz} and the entrance boundary condition \eqref{0.5'}, one can apply Theorem \ref{extension} and Corollary \ref{cor} to define processes solution to \eqref{sdeB} started from $\infty$. We sum up those results in the following corollary.
%Our last result provides sufficient conditions
\begin{corol}\label{NLCSBPstartedfrominfinity} Assume that Condition $\mathbb{H}$ holds. If \eqref{0.5'} and \eqref{onesidedLipschitz}  are satisfied, then  the flow of c\`adl\`ag strong solutions to \eqref{sdeB}, $(X^{(x)})_{x\in [0,\infty)}$,  converges weakly, in the Skorokhod topology, as $x\to\infty$, towards a non-explosive Feller process $X^{(\infty)}$ started from infinity.
%Moreover, for all $t\geq 0$, \[X_t^{(x)}\underset{x\rightarrow \infty}{\longrightarrow} X_t^{(\infty)}<\infty \text{ a.s.}\]
\end{corol}
\begin{remark} 
Condition $\mathbb{H}$ encompasses the comparison property 
\eqref{comparisonproperty}, so that under the assumptions of Corollary \ref{NLCSBPstartedfrominfinity}, for any $t>0$, $X_t^{(x)}$ increases almost-surely towards $X_t^{(\infty)}<\infty$, as $x$ goes to $\infty$.
\end{remark}
Sufficient conditions for continuous-state nonlinear branching processes to satisfy  \eqref{0.5'} are given in \cite[Theorem 2.13]{2017arXiv170801560L} and \cite[Theorem 2.2]{LePardoux19}. The problem of defining the process started from infinity was left unanswered in these works. Corollary \ref{NLCSBPstartedfrominfinity} enables us to address this question and to give a sense, in terms of sample paths, to the  phenomenon of coming down from infinity. 

To exemplify Corollary \ref{NLCSBPstartedfrominfinity}, consider the case for which $\gamma_i(x)=x$ for all $x\geq 0$, $i=1,2$ and $\gamma_0$ is continuous on $\mathbb{R}_+$. In this setting, the process solution to \eqref{sdeB} is also called continuous-state branching process with \textit{competition}. The nonlinear drift $\gamma_0$ is interpreted as modelling interactions between individuals, see Le and Pardoux \cite{LePardoux19}. The following corollary combines some results established in \cite{LePardoux19} and  Corollary \ref{NLCSBPstartedfrominfinity}.
\begin{corol}\label{CSBPwithcompetition} Assume that $\gamma_i(x)=x$ for all $x\geq 0$, $i=1,2$ and  that $\gamma_0$ is continuous on $\mathbb{R}_+$, with $\gamma_0(0)=0$. If $\gamma_0$ satisfies \eqref{onesidedLipschitz} and $\int^{\infty}_b\frac{\ddr x}{|\gamma_0(x)|}<\infty$ for some $b>0$, then Condition $\mathbb{H}$ and \eqref{0.5'} are satisfied, and the process $(X^{(\infty)}_t,t\geq 0)$, started from $\infty$, is well-defined, comes down from infinity instantaneously and then stays finite almost-surely. 
\end{corol} 
\begin{remark} \begin{enumerate}
\item An explicit example is the \textit{logistic} continuous-state branching process for which  for all $x\geq 0$, $\gamma_i(x)=x$, $i=1,2$ and $\gamma_0(x)=-\frac{c}{2}x^{2}$ with $c>0$. It satisfies \eqref{0.5'}, starts from infinity, ``comes down" instantaneously, and then stays finite almost-surely. We refer to \cite{Lambert}, see also \cite[Lemma 6.4]{Foucart} where \eqref{0.5'} is established.
\vspace{-2mm}
\item The non-explosion of the process in Corollary \ref{CSBPwithcompetition} will be guaranteed by the assumption $\int_0^{\infty}(z\wedge z^2)\nu(\ddr z)<\infty$. When the latter integral is infinite, different behaviors at $\infty$ may occur, see for instance \cite{Foucart} for the case $\gamma_0(x)=-\frac{c}{2}x^{2}$ for all $x\geq 0$.
\vspace{-2mm}
\item Examples of continuous-state nonlinear branching processes for which $\gamma_i(x)\neq x$, $i=1,2$ and Corollary \ref{NLCSBPstartedfrominfinity} applies can be found in \cite[Page 2536]{2017arXiv170801560L}. For instance, if $\nu(\ddr z)=c_{\alpha}\mathrm{1}_{\{z>0\}}z^{-1-\alpha} \ddr z$, with $\alpha\in (1,2)$, $c_{\alpha}>0$, $\gamma_0=\gamma_1=0$  and $\gamma_2(x)=x^{r_2}$ with $r_2>\alpha$, then the process solution to \eqref{sdeB} can be started from infinity and comes down from infinity. We refer to \cite{2016arXiv160909593L} for a different method allowing one to define this process started from infinity.
\end{enumerate}
\end{remark}
%\noindent As an explicit example, consider the process such that for all $x\geq 0$, $\gamma_i(x)=x$, $i=1,2$ and $\gamma_0(x)=-x^{r}$ with $r>1$. It satisfies \eqref{0.5'}, starts from infinity, ``comes down" instantaneously, and then stays finite almost-surely.}} 

\section{Proofs}
\subsection{Proof of Lemma \ref{equivcdi}}\label{prooflemmaequivcdi}
Recall the statement of Lemma \ref{equivcdi}. We show that $(a)\Longrightarrow (b)$. The condition \eqref{0.5'} entails that for any fixed $t$, there is $b$ large enough from which, $\underset{x\rightarrow \infty}{\limsup }\ \mathbb{P}_x(T_b> t)<1$. This implies that there exists $x_b\geq b$ such that $\underset{x\geq x_b}{\sup}\mathbb{P}_{x}(T_b> t)<1$. Since there are no negative jumps, for any $x\in [b, x_b]$, $\mathbb{P}_x(T_b>t)\leq \mathbb{P}_{x_b}(T_b>t)<1$. Therefore, if $b$ is large enough then, for all $x\geq b$,
\begin{equation}\label{boundentrance1} \mathbb{P}_x(T_b> t)\leq \sup_{x\geq b}\mathbb{P}_{x}(T_b> t)<1.
\end{equation}
%and any $x\geq x'_b$.
In particular, we see that $\mathbb{P}_x(T_b<\infty)>0$ for all $x\geq b$ and that  there is $t>0$ and a large enough $b$ such that $$\alpha_b:=\sup_{x\geq b}\mathbb{P}_x(T_b>t)<1.$$ We now work with those $t$ and $b$. Let $x>b$, and $n\in \mathbb{N}$, by using the Markov property at time $t$
\begin{align*}
\mathbb{P}_x(T_b>nt)&=\mathbb{E}_x\left[\mathbbm{1}_{\{T_b>t\}}\mathbbm{1}_{\{T_b\circ \theta_t>(n-1)t\}}\right]\\
&=\mathbb{E}_x\left[\mathbbm{1}_{\{T_b>t\}}\mathbb{E}_{X_t}\left[\mathbbm{1}_{\{T_b>(n-1)t\}}\right]\right]\\
&\leq \mathbb{P}_x(T_b>t) \ \underset{x\geq b}{\sup}\ \mathbb{P}_{x}(T_b>(n-1)t).
\end{align*}
Hence for all $n\geq 1$,
%\begin{align*}
$\underset{x\geq b}{\sup}\ \mathbb{P}_x(T_b>nt) \leq \alpha_b \ \underset{x\geq b}{\sup}\ \mathbb{P}_{x}(T_b>(n-1)t)$
%\end{align*}
and $\underset{x\geq b}{\sup}\ \mathbb{P}_x(T_b>nt) \leq \alpha_b^{n}.$
This entails $\sup_{x\geq b} \mathbb{E}_x(T_b)\leq \sum_{n=0}^{\infty}\alpha_b^{n}<\infty.$

We now show that $(b)\Longrightarrow (c)$. Note that the limits in $(c)$ are monotone increasing in $x$ and decreasing in $b$. According to $(b)$, for large enough $b$, $\mathbb{E}_x(T_b)<\infty$ for any $x\geq b$ so that $\mathbb{P}_x(T_b<\infty)=1$.  Let $x>x'>b$. By the absence of negative jumps and the strong Markov property
\[\mathbb{E}_{x}(T_b)=\mathbb{E}_{x}(T_{x'})+\mathbb{E}_{x'}(T_b).\]
Letting $x$  go towards $\infty$ and then $x'$ go towards $\infty$, we get
\[\underset{x\rightarrow \infty}{\lim} \mathbb{E}_{x}(T_b)=\underset{x'\rightarrow \infty}{\lim}\ \underset{x\rightarrow \infty}{\lim}\mathbb{E}_{x}(T_{x'})+\underset{x'\rightarrow \infty}{\lim}\mathbb{E}_{x'}(T_b).\]
By $(b)$, $\underset{x\rightarrow \infty}{\lim} \mathbb{E}_{x}(T_b)<\infty$ and thus
$\underset{x'\rightarrow \infty}{\lim}\ \underset{x\rightarrow \infty}{\lim}\mathbb{E}_{x}(T_{x'})=0.$ The last implication $(c)\Longrightarrow (d)$ is a simple consequence of the Markov inequality.\qed

%We will also establish an important regularity result for the first passage time $T_b$ and the convergence in law of the processes started from $x$ to that started from infinity, when $x$ goes to $\infty$, in the Skorokhod sense.

\subsection{Proof of Theorem \ref{extension}}
Recall the statement of Theorem \ref{extension} and its assumptions. The arguments to prove Theorem \ref{extension} are adapted from those of Kallenberg in \cite[Theorem 23.13]{Kallenberg}. Let $y>x>b$. Since the process $(X_t,t\geq 0)$ has no negative jumps, $T_b>T_x$ almost surely under $\mathbb{P}_y$. Moreover for any $t\geq 0$,  $\{T_b< t\}\subset \{T_x<t, T_b-T_x<t\}$  and
\begin{align*}
\mathbb{P}_y(T_b<t)&\leq \mathbb{P}_y(T_x\leq t, T_b-T_x<t)\\
&=\mathbb{P}_y(T_x< t)\mathbb{P}_x(T_b\leq t)\leq \mathbb{P}_x(T_b<t).
\end{align*}
This implies that $(\mathbb{P}_x(T_b<t),x\geq 0)$ admits a limit as $x$ goes to $\infty$ and that $(\mathbb{E}_x(T_b),x\geq b)$ is non-decreasing in $x$.  The inequality \eqref{boundentrance1} ensures that $\underset{x\rightarrow \infty}{\lim} \mathbb{P}_x(T_b<t)>0$. Recall also that Lemma \ref{equivcdi} states that (\ref{0.5'}) holds if and only if
%\begin{equation}\label{its1}\underset{b\rightarrow \infty}{\lim }\ \underset{x\rightarrow \infty}{\lim}\mathbb{P}_x(T_b<t)=1
%\end{equation}
%as well as
\begin{equation}\label{unifboundmoment}
\sup_{x\geq b}\mathbb{E}_x(T_b)\underset{b\rightarrow \infty}{\longrightarrow} 0.
\end{equation}
For any function $f\in C_b(E)$, one denotes by $||f||$ the supremum norm of $f$.  We now show that $(P_tf(x),x\geq 0)$ admits a limit as $x$ goes to $\infty$ for any $f\in C_b(E)$ with $E=[0,\infty)$. Fix $t\geq 0$, for any $x>b$
\begin{align*}
P_tf(x)&=\mathbb{E}_x[\mathbbm{1}_{\{T_b>t\}}f(X_t)]+\mathbb{E}_x[\mathbbm{1}_{\{T_b\leq t\}}f(X_t)]\\
&=\mathbb{E}_x[\mathbbm{1}_{\{T_b>t\}}f(X_t)]+\mathbb{E}_x[\mathbbm{1}_{\{T_b\leq t\}}P_{t-T_b}f(X_{T_b})]\\
&=\mathbb{E}_x[\mathbbm{1}_{\{T_b>t\}}f(X_t)]+\mathbb{E}_x[\mathbbm{1}_{\{T_b\leq t\}}P_{t-T_b}f(b)],
\end{align*}
where we have used the strong Markov property at $T_b\wedge t$ in the second equality and in the third, the absence of negative jumps which implies $X_{T_b}=b$ a.s.
Set \[g(a,s):=\mathbbm{1}_{s\leq t}P_{t-s}f(a)+\mathbbm{1}_{s>t}f(a).\] For any $x$ and $y$ larger than $b$,
\begin{align*}
|P_tf(y)-P_tf(x)|&=\lvert \mathbb{E}_x[\mathbbm{1}_{\{T_b>t\}}f(X_t)]-\mathbb{E}_y[\mathbbm{1}_{\{T_b>t\}}f(X_t)]\\
&\quad+\mathbb{E}_x[\mathbbm{1}_{\{T_b\leq t\}}P_{t-T_b}f(X_{T_b})]
-\mathbb{E}_y[\mathbbm{1}_{\{T_b\leq t\}}P_{t-T_b}f(X_{T_b})]\lvert\\
&\leq ||f||(\mathbb{P}_x(T_b>t)+\mathbb{P}_y(T_b>t))\\
&\qquad+\left\lvert \mathbb{E}_x\left[g(b,T_b)-f(b)\mathbbm{1}_{\{T_b>t\}}\right]-\mathbb{E}_y\left[g(b,T_b)-f(b)\mathbbm{1}_{\{T_b>t\}}\right]\right\lvert\\
&\leq 2||f||(\mathbb{P}_x(T_b>t)+\mathbb{P}_y(T_b>t))+ \left\lvert\mathbb{E}_x\left[g(b,T_b)\right]-\mathbb{E}_y\left[g(b,T_b)\right]\right\lvert.
\end{align*}From (\ref{0.5'}), we see that for $x,y\rightarrow \infty$ and then $b\rightarrow \infty$,
\[2||f||(\mathbb{P}_x(T_b>t)+\mathbb{P}_y(T_b>t))\longrightarrow 0.\]
  Moreover, (\ref{unifboundmoment}) entails $\sup_{x}\mathbb{E}_x(T_b)<\infty$
which provides  for any fixed $b$
$$\underset{s\rightarrow \infty}{\lim} \sup_{x\geq 0} \mathbb{P}_x(T_b>s)=0.$$
The family of laws of $T_b$ under $\mathbb{P}_x$ for $x\geq 0$, $(\mathbb{P}_x\circ T_b^{-1},x\geq 0)$, is therefore tight and admits a unique limit since $(\mathbb{P}_x(T_b<t),x\geq 0)$ converges. Thus,  $(\mathbb{P}_x\circ T_b^{-1},x\geq 0)$ converges weakly and since the map $s\mapsto g(b,s)$ is bounded and continuous, therefore, $$\left\lvert\mathbb{E}_x\left[g(b,T_b)\right]-\mathbb{E}_y\left[g(b,T_b)\right]\right\lvert\underset{x,y\rightarrow \infty}{\longrightarrow} 0.$$
Since $|P_tf(y)-P_tf(x)|\underset{x,y\rightarrow \infty}{\longrightarrow} 0$, then for any sequence $(x_n,n\geq 1)$ such that $x_n\underset{n\rightarrow \infty}{\longrightarrow} \infty$, $(P_tf(x_n),n\geq 1)$ is a Cauchy sequence and admits a limit in $\mathbb{R}$. This limit does not depend on the sequence $(x_n, n\geq 1)$ and we set
$P_tf(\infty):=\underset{x\rightarrow \infty}{\lim} P_tf(x)$. By the assumption of non-explosiveness, for all $x\in E$, the transition kernels $P_t(x,\cdot):=\mathbb{P}(X_t\in \cdot|X_0=x)$ are probability measures over $E=[0,\infty)$. Since the convergence $P_tf(\infty):=\underset{x\rightarrow \infty}{\lim} P_tf(x)$, holds for any $f\in C_b(E)$,  the probability measures $P_t(x,\cdot)$ over $E$ converge weakly towards $P_t(\infty,\cdot)$ as $x$ goes to $\infty$. Taking the constant function $f=1$, we see that  $P_t(\infty,\cdot)$ is a probability measure over $E$.

%We now extend $P_t$ on $(0,\infty]$ by letting $P_tf(\infty):=\underset{x\rightarrow \infty}{\lim} P_tf(x)$ for any $f\in C_b$ and $t>0$.
We proceed to check that $(P_t)$ forms a Feller semigroup on $\bar{E}=[0,\infty]$. Recall $C(\bar{E})$. Let $f\in C(\bar{E})$.  By the assumption $(ii)$,  $P_tf$ is continuous on $E$. By the definition, $P_tf(\infty)=\underset{x\rightarrow \infty}{\lim}P_tf(x)<\infty$ and $P_tf$ is continuous at $\infty$. Therefore, $P_t$ maps $C(\bar{E})$ into itself. By the assumption $\textit{(ii)}$, for any $x\in [0,\infty)$, one has $P_tf(x)\underset{t\rightarrow 0}{\longrightarrow} f(x)$, and it only remains to show that $P_tf(\infty)\underset{t\rightarrow 0}{\longrightarrow} f(\infty)$. Let $\varepsilon>0$ and choose $b$ large enough such that $\sup_{x\geq b}|f(x)-f(\infty)|\leq \varepsilon$. Note that there is $t$ small enough such that $\mathbb{P}_{2b}(T_b\leq t)\leq \varepsilon$. Then,
\begin{align*}
|P_tf(\infty)-f(\infty)|&\leq \underset{x\rightarrow \infty}\lim \big(\mathbb{E}_x[|f(X_t)-f(\infty)|\mathbbm{1}_{\{T_b\leq t\}}]+\mathbb{E}_x[|f(X_t)-f(\infty)|\mathbbm{1}_{\{T_b>t\}}]\big)\\
&\leq 2||f||\mathbb{P}_{2b}(T_b\leq t)+\varepsilon\\
&\leq (2||f||+1)\varepsilon.
\end{align*}
The pointwise continuity of the semigroup at $t=0$ is therefore satisfied on $[0,\infty]$ and finally the semigroup on $\bar{E}$ satisfies $\textit{(i)}$ and $\textit{(ii)}$ and is Feller. Theorem 19.15 in \cite{Kallenberg} provides the existence of a strong Markov process $(X_t,t\geq 0)$ with c\`adl\`ag paths started from $\infty$. We denote its law by $\mathbb{P}_\infty$. Note that for any $t>0$, $P_t(\infty, E)=\mathbb{P}_\infty(X_t<\infty)=1$. It remains to show that $\mathbb{P}_\infty(\forall t>0, X_t<\infty)=1$. By the assumption for any $x\in E$, $(X_t,t\geq 0)$ does not explode under $\mathbb{P}_x$. Let $s>0$. By the Markov property at time $s$ under $\mathbb{P}_\infty$, \[\mathbb{P}_\infty(\exists t>s; X_{t}=\infty)=\mathbb{E}_{\infty}(\mathbb{P}_{X_s}(\exists t>0; \tilde{X}_t=\infty))\] where $\tilde{X}$ is a copy of $X$. Since $\mathbb{P}_\infty(X_s<\infty)=1$, we have that for any $s>0$,  $\mathbb{P}_\infty$-a.s. $\mathbb{P}_{X_s}(\exists t>0; \tilde{X}_t=\infty)=0$. Hence $\mathbb{P}_\infty(\exists t>s; X_{t}=\infty)=0$, which allows one to conclude since $s$ is arbitrary. \qed

\subsection{Proof of Proposition \ref{momentunderPinfty}}\label{A2}
Assume that \eqref{0.5'} holds and recall $\mathbb{E}_\infty$ the corresponding expectation under $\mathbb{P}_\infty$. We show first the exponential moment property stated in (a). Recall that under $\mathbb{P}_\infty$, $T_b$ decreases as $b$ goes to $\infty$. Denote by $T_\infty=\inf \{t>0; X_t<\infty\}$ its limit. By Theorem \ref{extension},
%, for any $t>0$
%\[\mathbb{P}_\infty(X_t<\infty)=\mathbb{P}_\infty(T_\infty\leq t)=1.\]
%Therefore, 
$T_\infty=0$ $\mathbb{P}_\infty$-a.s. Let $\theta>0$. Fix $t>0$. By Lebesgue's theorem
\beqlb\label{1}
\lim_{b\to\infty}\mbb{P}_\infty(T_b>t)=0.
\eeqlb
Hence, there exists $b_\theta>0$, such that for all $b\geq b_\theta$, $e^{\theta t}\mbb{P}_\infty(T_b> t)<1$. Moreover for any $n\geq 1$, by the Markov property, we get $\mbb{P}_\infty(T_b> nt)\leq \mbb{P}_\infty(T_b> t)^n$ and
%and
%\beqnn
%\mbb{P}_\infty(T_b> nt)
%\ar=\ar
%\mbb{E}_\infty\Big[\mathbbm{1}_{\{T_b> t\}}\mathbbm{1}_{\{T_b\circ\theta_t> (n-1)t\}}\Big] \cr
% \ar=\ar
%\mbb{E}_\infty\Big[\mathbbm{1}_{\{T_b> t\}}\mbb{E}_{X_t}[\mathbbm{1}_{\{T_b> (n-1)t\}}]\Big] \cr
% \ar\le\ar
%\mbb{E}_\infty \Big[\mathbbm{1}_{\{T_b> t\}}\mbb{E}_\infty [\mathbbm{1}_{\{T_b> (n-1)t\}}]\Big] \cr
% \ar=\ar
%\mbb{P}_\infty(T_b> t)\mbb{P}_\infty(T_b> (n-1)t)\cr
% \ar\le\ar
%\mbb{P}_\infty(T_b> t)^n.
% \eeqnn
for any $b\geq b_\theta$
\begin{align*}
\mathbb{E}_\infty(e^{\theta T_b})= \mathbb{E}_\infty\left[\theta\int_{0}^{T_b}e^{\theta s}\ddr s\right]+1&=\theta t \int_{0}^{\infty}e^{\theta t s}\mathbb{P}_\infty(T_b>ts)\ddr s+1\\
&=\theta t\sum_{n=0}^{\infty}\int_{n}^{n+1}e^{\theta ts} \mathbb{P}_\infty(T_b>ts)\ddr s+1\\
&\leq \theta te^{\theta t}\sum_{n=0}^{\infty}\left( e^{\theta t}\mathbb{P}_\infty(T_b>t)\right)^{n}+1<\infty.
\end{align*}
%\[
%\mathbb{E}_\infty(T_b)=t\int_{0}^{\infty}\mathbb{P}_\infty(T_b>ts)\ddr s=t\sum_{n=0}^{\infty}\int_{n}^{n+1}\mathbb{P}_\infty(T_b>ts)\ddr s \leq t\sum_{n=0}^{\infty}\mathbb{P}_\infty(T_b>nt)<\infty.
%\]
%\begin{align*}
%\mathbb{E}_\infty(T_b)
%%&=t\int_{0}^{\infty}\mathbb{P}_\infty(T_b>ts)\ddr s\\
%%&=t\sum_{n=0}^{\infty}\int_{n}^{n+1}\mathbb{P}_\infty(T_b>ts)\ddr s\\
%&\leq t\sum_{n=0}^{\infty}\mathbb{P}_\infty(T_b>nt)<\infty.
%\end{align*}
We establish now the convergence in \eqref{convergencemoment}.
Denote by $\theta$ the shift operator, see e.g. \cite[p.146]{Kallenberg}. For any $x>b>0$,
\begin{equation}\label{ap1}
T_b\circ\theta_{T_x}=T_b-T_x\longrightarrow T_b\quad\mbox{as}\quad x\to\infty,\ \mbb{P}_\infty\mbox{ -a.s.}
\end{equation}
By the strong Markov property,
\begin{equation*}
\mbb{E}_\infty [h(T_b)]=\mbb{E}_\infty[h(T_b)-h(T_b\circ\theta_{T_x})]+\mbb{E}_x[(h(T_b))].
\end{equation*}
Recall \eqref{ap1} and that $h$ is either bounded or  nonnegative and increasing. By
applying Lebesgue's theorem, we get
\begin{equation*}
  \mbb{E}_\infty [h(T_b)]=\lim_{x\to\infty}\mbb{E}_x[h(T_b)].
\end{equation*}
 \qed

\subsection{Proof of Theorem \ref{weakconv}}
Define the metric $\rho$ on $\bar{E}=[0,\infty]$ by $\rho(x,y)=|e^{-x}-e^{-y}|$ for any $x,y\in \bar{E}$ and let $D$ be the space of c\`adl\`ag functions $f: E\to \bar{E}$. We endow $D$ with the Skorokhod topology, for which we refer, for instance, to Ethier and Kurtz's book \cite[Chapter 3, p.116]{EthierKurtz}. The proof of Theorem \ref{weakconv} follows by combining Lemma \ref{tightness} (tightness) and Lemma \ref{finite} (convergence of finite dimensional laws), given below, and by applying \cite[Corollary 9.3]{EthierKurtz}.

%Recall that we denote by $C([0,\infty])$ the space of continuous functions on $[0,\infty)$ with a finite limit at $\infty$.
Let $(P_t)_{t\ge0}$ be a Feller semigroup on $C(\bar{E})$.
%Note that   on $C_b(\bar{\mbb{R}}_+)=C_0(\bar{\mbb{R}}_+)$.
Let $A:\mcr{D}\to C(\bar{E})$ be the generator of $(P_t)_{t\ge0}$ with $\mcr{D}$ its domain. Then by \cite[Theorem 19.4]{Kallenberg}, we see that  $\mcr{D}$ is dense in $C(\bar{E})$ with respect to the supremum norm $||\cdot||$. For each $x\in\bar{E}$, denote by $X^{(x)}=(X^{(x)}_t)_{t\ge 0}$ the corresponding Feller process  with initial value $x>0$ and càdlàg paths. From Dynkin's formula, see \cite[Theorem 19.21]{Kallenberg}, we have for any bounded stopping time $\tau$, and any $f\in\mcr{D}$
\begin{equation*}
\mbb{E}[f(X^{(x)}_\tau)]=f(x)+\mbb{E}\left[\int^\tau_0Af(X^{(x)}_s)ds\right].
\end{equation*}

\blemma\label{tightness} The sequence of processes $(X^{(x)})_{x>0}$ is relatively compact in distribution, i.e. every subsequence has a subsequence that converges in distribution in $D$.
\elemma
\proof
Since $\mcr{D}$ is dense in $C(\bar{E})$, for any $\varepsilon>0$ and $f\in\mcr{D}$, there exists $g\in\mcr{D}$ such that $||f^2-g||<\varepsilon$. Applying the Markov property and Dynkin's formula we have  for any
bounded stopping time $\tau$ and constant $h>0$,
\begin{align*}
\mathbb{E}[(f(X^{(x)}_{\tau+h})-f(X^{(x)}_{\tau}))^2]
&=\mathbb{E}[f^2(X^{(x)}_{\tau+h})-f^2(X^{(x)}_{\tau})]\\
&
\qquad -
2\mathbb{E}\left[f(X^{(x)}_\tau)\mbb{E}[f(X^{(x)}_{\tau+h})-f(X^{(x)}_\tau)|\mcr{F}_\tau]\right]\\
&\le 2\varepsilon+\mathbb{E}[g(X^{(x)}_{\tau+h})-g(X^{(x)}_{\tau})]+2||f|| \left \lvert \mbb{E}\left(\mbb{E}[f(X^{(x)}_{\tau+h})-f(X^{(x)}_\tau)\lvert \mcr{F}_{\tau}]\right) \right \lvert\\
&= 2\varepsilon +\mathbb{E}\left[\int^{\tau+h}_\tau Ag(X^{(x)}_s)ds\right]+2||f|| \left \lvert \mbb{E}\left[\int^{\tau+h}_\tau Af(X^{(x)}_s)ds\right]\right \lvert\\
&\le
2\varepsilon+||Ag||h+2||f||\cdot||Af||h.
\end{align*}
Therefore, by Jensen's inequality and the arbitrariness of $\varepsilon$,
$$
\lim_{\delta\to0}\limsup_{x\to\infty}\sup_{\tau\le t}\sup_{h\in[0,\delta]}\mbb{E}|f(X^{(x)}_{\tau+h})-f(X^{(x)}_{\tau})|=0\qquad t>0,
$$
where the second supremum is over all stopping time $\tau\leq t$ for the process $X^{(x)}$. Then by Aldous' criterion of tightness, see e.g. \cite[Theorem 16.11 and Lemma 16.12, p.314]{Kallenberg}, for any $f\in\mcr{D}$, $(f\circ X^{(x)})_{x>0}$ is  relatively compact in distribution. Finally, since $\mcr{D}$ is dense in $C(\bar{E})$, applying \cite[Theorem 9.1 p.142]{EthierKurtz} we see that  $(X^{(x)})_{x>0}$ is relatively compact in distribution. Note that $\bar{E}$ is compact, thus the compact containment condition, needed for \cite[Theorem 9.1 p.142]{EthierKurtz} to apply, holds. \qed

\blemma\label{finite}
The finite-dimensional distributions of $X^{(x)}$ converges weakly to that of $X^{(\infty)}$ as $x\to\infty.$
\elemma
\proof
We are going to prove that for any sequence $\{t_1<t_2<\dots<t_n\}\subset E$ and $\{f_1,f_2,\dots,f_n\}\subset C({\bar{E}})$,
\begin{equation}\label{fd1}
\lim_{x\to\infty}\mbb{E}[f_1(X^{(x)}_{t_1})f_2(X^{(x)}_{t_2})\dots f_n(X^{(x)}_{t_n})]=\mbb{E}[f_1(X^{(\infty)}_{t_1})f_2(X^{(\infty)}_{t_2})\dots f_n(X^{(\infty)}_{t_n})].
\end{equation}
If (\ref{fd1}) holds, then by the Stone-Weierstrass theorem, see e.g. \cite[Theorem 7.29]{HS} for any $F_n\in C({\bar{E}^n})$ we have
$$
\lim_{x\to\infty}\mbb{E}[F_n(X^{(x)}_{t_1},X^{(x)}_{t_2},\dots ,X^{(x)}_{t_n})]=\mbb{E}[F_n(X^{(\infty)}_{t_1},X^{(\infty)}_{t_2},\dots, X^{(\infty)}_{t_n})].
$$
Obviously, (\ref{fd1}) holds for $n=1$. Assume that (\ref{fd1}) holds for $n=m$. Then by the Markov property
\begin{align*}
&\mbb{E}[f_1(X^{(x)}_{t_1})f_2(X^{(x)}_{t_2})\dots f_{m+1}(X^{(x)}_{t_{m+1}})]\\
& =\mbb{E}[f_1(X^{(x)}_{t_1})f_2(X^{(x)}_{t_2})\dots f_m(X^{(x)}_{t_m})\mbb{E}[f_{m+1}(X^{(x)}_{t_{m+1}})|\mcr{F}_{t_{m}}]]\\
& =\mbb{E}[f_1(X^{(x)}_{t_1})f_2(X^{(x)}_{t_2})\dots f_m(X^{(x)}_{t_m})P_{t_{m+1}-t_m}f_{t_{m+1}}(X^{x}_{t_m})].
\end{align*}\\
By the Feller property, the map $x\mapsto f_m(x)P_{t_{m+1}-t_m}f_{t_{m+1}}(x)$ is continuous and bounded.
Then we can finish the proof of (\ref{fd1}) by induction.\qed
\subsection{Proof of Proposition \ref{Fellerprop}}
%Let $\mcr{P}$ be the set of probability measures $\mu$ such that $\int_E z\mu(\d z)<\infty$.
%For any two probability measures $\mu_1, \mu_2\in\mcr{P}$, the $L^1$-Wasserstein distance $W_1$ is defined by $$
%W_{1}(\mu_1,\mu_2)=\inf_{\Pi\in\mathscr{C}(\mu_1,
%	\mu_2)}\int_{E^2} |x-y|\,\Pi(\d x, \d y), $$ where
%$\mathscr{C}(\mu_1, \mu_2)$ is the collection of measures on
%$E^2$ with marginals $\mu_1$ and $\mu_2$.
%We denote by $(X_t(x))_{t\ge0}$ the strong solution to (\ref{sdeB}) with initial value $x$. Using order persevering property, we see
%  $X_t(x)\leq X_t(y)$ a.s. for any $t\geq 0$ and $0\le x<y$. Then for any $0\le x<y$,
%\begin{align*}
%W_1(P_t(x,\d z), P_t(y, \d z))&\le
%\mathbb{E}(|X_t(y)-X_t(x)|)\\
%&=\mathbb{E}(X_t(y)-X_t(x))\\
%&=y-x+\int_{0}^{t}\mathbb{E}(\gamma_0(X_s(y))-\gamma_0(X_s(x)))\ddr s\\
%&\leq y-x+\theta \int_{0}^{t}\mathbb{E}(X_s(y)-X_s(x))\ddr s.
%\end{align*}
%By Gronwall's Lemma,
%$$W_1(P_t(x,\d z), P_t(y, \d z))\leq e^{\theta t}|x-y|$$ for all $t$, and we can use \cite[Lemma 5.3]{Chen04} to establish condition \textit {(ii)}.
Recall $(X_t^{(x)}, t\geq 0)$ the solution to \eqref{sdeB} and the statement of Proposition \ref{Fellerprop}. Let $b\mcr{L}$ be the set of bounded Lipschitz continuous functions on $E$ and $L(f)$ denotes the Lipschitz constant of $f\in b\mcr{L}$.
For $f\in C_b(E)$ and $M>0$, let
\begin{equation*}
f_M(x)=\begin{cases}f(x),\quad &0\le x\le M,\\
f(M),\quad & x>M.
\end{cases}
\end{equation*}
It follows from the Stone-Weierstrass theorem that there exists a sequence of bounded Lipschitz functions $(f^{(n)}_M)_{n\geq 1}$ such that
\begin{equation}\label{converge1}
||f^{(n)}_M-f_M||\to0\qquad\mbox{as}~n\to\infty.
\end{equation}
Since by the assumption, $X_t^{(x)}\leq X_t^{(y)}$ a.s. for any $t\geq 0$ and $0\le x<y$,  for any $0\le x<y$, by \eqref{onesidedLipschitz}
\begin{align*}
\mathbb{E}(|X_t^{(y)}-X_t^{(x)}|)=\mathbb{E}(X_t^{(y)}-X_t^{(x)})
&=y-x+\int_{0}^{t}\mathbb{E}(\gamma_0(X_s^{(y)})-\gamma_0(X_s^{(x)}))\ddr s\\
&\leq y-x+\theta \int_{0}^{t}\mathbb{E}(X_s^{(y)}-X_s^{(x)})\ddr s.\end{align*}
By Gronwall's Lemma,
$$\mathbb{E}(|X_t^{(y)}-X_t^{(x)}|) \leq e^{\theta t}|x-y|$$ for all $t$. Thus, for any $f\in b\mcr{L}$,
\beqlb\label{Lip}
|P_tf(x)-P_tf(y)|\le L(f)e^{\theta t} |x-y|,
\eeqlb
which entails that $P_tf\in b\mcr{L}$. Define
\begin{equation*}
g^   M(x)=\begin{cases}
0,\quad &0\le x\le M/2,\\
\frac{2x}{M}-1,\quad &M/2\le x\le M,\\
1,\quad & x>M.
\end{cases}
\end{equation*}
Fix $t$ and $x$ and denote by $P_t(x, \cdot)$ the law of $X^{(x)}_t$. For any $\varepsilon>0$, we can choose $M, \delta>0$, such that for each $y\in((x-\delta)\vee0, x+\delta)$
\begin{align}\label{es1}
|P_tf(y)-P_tf_M(y)|&\le\nonumber
2||f|| P_t(y, [M, \infty))\\\nonumber
&\le
2||f|| \left(P_tg^M(x)+|P_t g^M(y)-P_t g^M(x)|\right)\\
&\le
2||f||\Big(P_t(x, [M/2, \infty))+\frac{2e^{\theta t}}{M}|x-y|\Big)
<\varepsilon/4
\end{align}
where, for any Borel set $A$, $P_t(y,A):=P_t\mathbbm{1}_{A}(y)$ and for the third inequality we used (\ref{Lip}).
From (\ref{converge1}), we see that there exists $n\in \mathbb{N}$ such that
\beqlb\label{es2}
||f_M-f_M^{(n)}||\le\varepsilon/8.
\eeqlb
Since $P_tf_M^{(n)}\in b\mcr{L}$, there is $\delta'<\delta$ such that for $y\in ((x-\delta')\vee 0, y+\delta'),$
\beqlb\label{es3}
\big|P_tf_M^{(n)}(x)-P_tf_M^{(n)}(y)\big|\le \varepsilon/4
\eeqlb
Combining (\ref{es1}), (\ref{es2}) and (\ref{es3}), we see that for $y\in ((x-\delta')\vee 0, y+\delta'),$
\beqnn
|P_tf(x)-P_tf(y)|\ar\le\ar|P_tf_M(x)-P_tf_M(y)|+|P_tf(x)-P_tf_M(x)|+|P_tf(y)-P_tf_M(y)|\\
\ar\le\ar
\big|P_tf_M^{(n)}(x)-P_tf_M^{(n)}(y)\big|+2||f_M-f^{(n)}_M||\\
\ar\ar\quad+|P_tf(x)-P_tf_M(x)|+|P_tf(y)-P_tf_M(y)|\\
\ar\le\ar
\varepsilon/4+\varepsilon/4+\varepsilon/4+\varepsilon/4=\varepsilon.
\eeqnn
Thanks to the arbitrariness of $\varepsilon$, we can complete the proof.\qed
\subsection{Proof of Corollaries \ref{NLCSBPstartedfrominfinity} an \ref{CSBPwithcompetition}}
\subsubsection{Proof of Corollary \ref{NLCSBPstartedfrominfinity}}
Recall Corollary \ref{NLCSBPstartedfrominfinity}. Assume that \eqref{0.5'} is satisfied. Under Condition $\mathbb{H}$, the solution of \eqref{sdeB} has no negative jump and is right-continuous, so that condition \textit{(ii)} is satisfied. By Proposition \ref{Fellerprop}, under the assumption \eqref{onesidedLipschitz}, the process $(X_t,t\geq 0)$ satisfies \textit{(i)}. Assumptions of Theorem \ref{extension} are therefore fulfilled and the process admits a Feller extension started from infinity. The weak convergence is a consequence of Corollary \ref{cor}.
%Moreover, by \eqref{comparisonproperty}, for any $t\geq 0$, $X_t^{(x)}$ increases almost-surely towards $X_t^{(\infty)}$ as $x$ goes to $\infty$.
\qed
\subsubsection{Proof of Corollary \ref{CSBPwithcompetition}}
Consider now the continuous-state branching processes with competition. These processes are solution to the SDE \eqref{sdeB} with $\gamma_i(x)=x$ for all $x\geq 0$ and $i=1,2$, with a continuous function $\gamma_0$ satisfying $\gamma_0(0)=0$.

Assume that $\gamma_0$ satisfies the one-sided Lipschitz condition \eqref{onesidedLipschitz}. We establish that \eqref{onesidedLipschitz} combined with $\int^{\infty}_0 (z\wedge z^{2})\nu(\ddr z)<\infty$, entails Condition $\mathbb{H}$. Dawson and Li \cite[Theorem 2.1]{DL} ensures that Equation \eqref{sdeB} admits a unique nonnegative c\`adl\`ag solution $(X_t^{(x)}, t\geq 0)$, started from $x\in [0,\infty)$. Moreover, the assumption $\gamma_0(0)=0$, combined with \eqref{onesidedLipschitz} entails that $\gamma_0(x)\leq \theta x$ for all $x\geq 0$. Assumptions of \cite[Theorem 2.2]{DL} are thus satisfied and the comparison property \eqref{comparisonproperty} holds true. A second application of \cite[Theorem 2.2]{DL} ensures that for all $x\geq 0$ and all $t\geq 0$, $X_t^{(x)}\leq Z^{(x)}_t$ where $(Z^{(x)}_t,t\geq 0)$ is a supercritical continuous-state branching process with L\'evy measure $\nu$, started from $x$. By the assumption $\int^{\infty}_0 (z\wedge z^{2})\nu(\ddr z)<\infty$, the process $(Z^{(x)}_t, t\geq 0)$ has finite mean and thus is non-explosive. By comparison, the process $(X^{x)}_t, t\geq 0)$ is thus also non-explosive for any $x\in \mathbb{R}_+$, and finally Condition $\mathbb{H}$ holds true. 

It remains to establish \eqref{0.5'}. Note that \eqref{onesidedLipschitz} is hypothesis (H1) in \cite{LePardoux19}. Assume that  $\gamma_0$ satisfies $\gamma_0(0)=0$ and \eqref{onesidedLipschitz}. By \cite[Theorem 2.2]{LePardoux19} and its proof, see \cite[page 10]{LePardoux19}, for any $b>0$,
%\begin{align*}
if $\int^{\infty}_b\frac{\ddr x}{|\gamma_0(x)|}<\infty$, then $\underset{x\geq b}\sup \ \mathbb{E}_x(T_b)<\infty$. Thus, by Lemma \ref{equivcdi}, if $\int^{\infty}_b\frac{\ddr x}{|\gamma_0(x)|}<\infty$ for some $b>0$, then \eqref{0.5'} is satisfied.  Hence Theorem \ref{extension} and Corollary \ref{NLCSBPstartedfrominfinity} apply.\qed 

\noindent \textbf{Acknowledgements:} C.F's research is partially  supported by LABEX MME-DII (ANR11-LBX-0023-01). P.L. and X.Z.'s research is supported by Natural Sciences and Engineering Research Council of Canada (RGPIN-2016-06704).
P.L.'s research is supported by Natural Sciences and Engineering Research Council of Canada (RGPIN-2012-07750) and National Natural Science Foundation of China (No.\ 11901570 and No.\ 11771046). X.Z.'s research is supported by  National Natural Science Foundation of China (No.\  11731012). The authors are grateful to Donald Dawson for his support and encouragement. C.F would like to thank Leif D\"oring for mentioning the proof of Kallenberg during the 3th Workshop on Branching Processes and Related Topics held in 2017 in Beijing.
\bibliographystyle{plain}
%\bibliography{doku}

\end{document}